\documentclass{amsart}

\usepackage[mathscr]{eucal}
\usepackage{amsmath}
\usepackage{amssymb}
\usepackage{amscd}
\usepackage{amsfonts}

\newtheorem{thm}{Donotwrite}[section]

\newtheorem{definition}[thm]{Definition}
\newtheorem{theorem}[thm]{Theorem}
\newtheorem{proposition}[thm]{Proposition}
\newtheorem{lemma}[thm]{Lemma}

\numberwithin{equation}{section}

\newfont{\germ}{eufm10}

\newcommand\cb{{\bf c}}
\newcommand\cd{\cdots}
\newcommand\eps{\epsilon}
\newcommand\et[1]{\tilde{e}_{#1}}
\newcommand\ft[1]{\tilde{f}_{#1}}
\newcommand\geh{\goth{g}}

\newcommand\goth[1]{\mbox{\germ #1}}
\newcommand\Kz{K_{\Z}}
\newcommand\La{\Lambda}
\newcommand\la{\lambda}
\newcommand\ol{\overline}
\newcommand\olLa{\ol{\La}}
\newcommand\ot{\otimes}
\newcommand\Q{{\mathbb Q}}
\newcommand\Uq{U_q(\geh)}
\newcommand\Uqp{U'_q(\geh)}
\newcommand\veps{\varepsilon}
\newcommand\vphi{\varphi}
\newcommand\wt{\mbox{\sl wt}\,}
\newcommand\Z{{\mathbb Z}}
\newcommand\Zn{\Z_{\ge0}}

\title[Existence of Crystal Bases]{Existence of Crystal Bases for \\
 Kirillov-Reshetikhin Modules of Type $D$} 

\author{Masato Okado}
\address{
Department of Mathematical Science, Graduate School of Engineering Science, 
Osaka University, Toyonaka, Osaka 560-8531, Japan}

\date{}

\begin{document}

\maketitle


\section{Introduction}

Let $\geh$ be an affine algebra and let $\Uqp$ be the corresponding quantum affine algebra without
degree operator. Among irreducible finite-dimensional $\Uqp$-modules there exists a distinguished
family called Kirillov-Reshetikhin modules (KR modules for short). They were introduced in \cite{KR}
in connection with a certain conjectural formula of multiplicities of irreducible 
$U_q(\geh_0)$-modules in a tensor product of those modules. Here $\geh_0$ stands for the 
finite-dimensional simple Lie algebra whose Dynkin diagram is obtained by removing the $0$-vertex,
that is prescribed in \cite{Kac}, from that of $\geh$. It is known \cite{D,CP} that irreducible 
finite-dimensional $\Uqp$-modules are classified by $n$-tuples of polynomials called Drinfeld 
polynomials, where $n$ is the rank of $\geh_0$.

Let us define KR modules by the Drinfeld polynomials. Let $k\in\{1,2,\ldots,n\},l\in\Z_{>0}$ and
$a$ an invertible element of $\Q(q)$. A KR module $\tilde{W}_{l,a}^{(k)}$ is defined to be the 
unique irreducible finite-dimensional $\Uqp$-module that has
\[
P_i(u)=\left\{
\begin{array}{ll}
(1-aq_i^{1-l}u)(1-aq_i^{3-l}u)\cd(1-aq_i^{l-1}u)\quad&\mbox{if }i=k,\\
1&\mbox{otherwise}
\end{array}\right.
\]
as Drinfeld polynomials. (See section \ref{subsec:crystal base} for the definition of $q_i$.)
When $l=1$ it is also called a fundamental representation. Since a fundamental representation
is known to have a crystal base \cite{K3}, we choose $a^\dagger$ so that 
${\tilde W}_{l,a^\dagger}^{(k)}$ has a crystal base and redefine $W_{l,a}^{(k)}=
{\tilde W}_{l,a^\dagger a}^{(k)}$. $W_{l,1}^{(k)}$ is also denoted simply by $W_l^{(k)}$.

The above-mentioned conjectural multiplicity formula was shown when $\geh$ is non-twisted by 
combining two results. The first one is the proof of certain algebraic relations among characters
of KR modules, called $Q$-systems, presented in \cite{N} for simply-laced cases and in \cite{H} 
for all non-twisted cases. The second one is a derivation of the multiplicity formula, called 
fermionic formula, in \cite{Ki} for type $A$ and in \cite{HKOTY} for all non-twisted cases,
by using the $Q$-systems.
However, it deserves to emphasize that there is also a $q$-analogue of the conjecture, called
$X=M$ conjecture \cite{HKOTY,HKOTT}. The definition of $X$ requires the existence of the crystal
base of a KR module. Despite many efforts as in \cite{KMN2,KKM,Y,Ko,JMO,K3,NS,SS,BFKL}, this 
existence problem is yet to be settled. For type $D$ for instance, the crystal base has been 
shown to exist for $W_l^{(k)}$ where $k=1,n-1,n;l\in\Z_{>0}$ in \cite{KMN2} and for $W_1^{(k)}$ for
arbitrary $k$ in \cite{K3}.

In this paper we prove the following theorem, thereby settling the problem for type $D$.

\begin{theorem} \label{th:main}
For $2\le k\le n-2$ and $l\ge1$, the $U'_q(D_n^{(1)})$-module $W_l^{(k)}$ has a crystal pseudobase.
\end{theorem}
\noindent
Here $(L,B)$ is said to be a crystal pseudobase if $(L,B/\{\pm1\})$ is a crystal base. 
(See Definition \ref{def:crystal base} for the definition of a crystal base.) Let us give 
a short sketch of our proof. We follow the technique already developed in \cite{KMN2}, namely,
from a fundamental representation $W_1^{(k)}$ we construct $W_l^{(k)}$ for any $l$ by fusion 
construction (section \ref{subsec:fusion}), and we apply a criterion of the existence of a crystal
pseudobase (Proposition \ref{prop:key}) to the constructed module $W_l^{(k)}$. Practically, we need
to check two conditions (\eqref{char} and \eqref{norm}). Checking the second one is not difficult, 
if once we find out the vectors $\{u_j\}$ correctly, whereas
checking the first one requires information on the image $W$ and the kernel $N$ of the $R$-matrix
$R(x,y):W_{1,x}^{(k)}\ot W_{1,y}^{(k)}\longrightarrow W_{1,y}^{(k)}\ot W_{1,x}^{(k)}$ at $x/y=q^2$.
Up to now such information was obtained by calculating the spectral decomposition of the $R$-matrix
when dealing with $W_l^{(k)}$ for higher $l$. It seems to be the reason why showing the existence of 
crystal bases of $W_l^{(k)}$ for higher $k$ has not been succeeded so far, since the calculation of
the corresponding $R$-matrix is too much complicated. However, thanks to the result by Nakajima 
\cite{N}, we are now able to identify $W$ and $N$ with tensor products of KR modules (Lemma 
\ref{lem:WN}). Using the crystal base of $W$ and a property of a bilinear form between 
$W_{1,q}^{(k)}\ot W_{1,q^{-1}}^{(k)}$ and $W_{1,q^{-1}}^{(k)}\ot W_{1,q}^{(k)}$, we can check the 
first condition of the criterion. It is known that a $\Uq$-module with a connected crystal base is 
irreducible. Therefore, once $W_l^{(k)}$ is shown to have a crystal pseudobase, it follows that 
$W_l^{(k)}$ is an irreducible finite-dimensional module with the desired Drinfeld polynomials, since
it is a simple quotient of $W_{1,q^{1-l}}^{(k)}\ot W_{1,q^{3-l}}^{(k)}\ot\cd\ot W_{1,q^{l-1}}^{(k)}$.

After the author finished the manuscript, he learned from Kashiwara that the module $W_l^{(k)}$ can
be shown to be irreducible by Theorem 9.2 of \cite{K3}. Once it is established, the character is 
known by \cite{C}. Hence it turns out that there is no need to prove the inequality of the character
in (i) just after Proposition \ref{prop:norm}. However, this does not seem to prove that $W_l^{(k)}$
is isomorphic to a module of the form of $V^{\ot l}/\sum_{i=0}^{l-2}V^{\ot i}\ot N\ot V^{\ot(l-2-i)}$.
The author was also informed from Nakajima that the existence of a polarization on the fundamental 
representation $W_1^{(k)}$ was shown in \cite{VV} (see also \cite{N2,BN} for more general results).
Hence similar calculations of the prepolarization as in section \ref{sec:pr2} will give a proof of 
the existence of crystal bases for KR modules of other quantum affine algebras.

\section{Crystal base and fusion construction}

\subsection{Crystal base} \label{subsec:crystal base}

In this subsection we briefly recall the definition of crystal bases. For more details along
with the definition of $\Uq$, refer to \cite{K1}.

Let $\geh$ be a symmetrizable Kac-Moody Lie algebra and let $M$ be a $\Uq$-module.
$M$ is said to be integrable if $M=\oplus_{\la\in P}M_{\la}$, $\dim M_{\la}<\infty$ for any $\la$, 
and for any $i$, $M$ is a union of finite-dimensional $U_q(\geh_i)$-modules.
Here $P$ is the weight lattice of $\geh$, $M_\la$ is the weight space of $M$ of
weight $\la$ and $U_q(\geh_i)$ is the subalgebra generated by Chevalley generators 
$e_i$ and $f_i$. If $M$ is integrable, we have
\begin{equation} \label{sl2decomp1}
M=\bigoplus_{0\le n\le\langle h_i,\la\rangle} f_i^{(n)}(\text{Ker}\;e_i\cap M_{\la}).
\end{equation}
Note that we use the following notations: $[m]_i=(q_i^m-q_i^{-m})/(q_i-q_i^{-1}),
[n]_i!=\prod_{m=1}^n[m]_i,f_i^{(n)}=f_i^n/[n]_i!$ with $q_i=q^{(\alpha_i,\alpha_i)}$,
where $(\;,\;)$ is an invariant bilinear form on $P$.
We define the endomorphisms $\et{i},\ft{i}$ of $M$ by
\begin{equation} \label{Kop1}
\ft{i}(f_i^{(n)}u)=f_i^{(n+1)}u\quad\text{and}\quad \et{i}(f_i^{(n)}u)=f_i^{(n-1)}u
\end{equation}
for $u\in\text{Ker}\;e_i\cap M_{\la}$ with $0\le n\le\langle h_i,\la\rangle$. Similarly, we have
\begin{equation} \label{sl2decomp2}
M=\bigoplus_{0\le n\le-\langle h_i,\mu\rangle} e_i^{(n)}(\text{Ker}\;f_i\cap M_{\mu}).
\end{equation}
These two decompositions are related as follows:
\begin{align*}
&\text{if }0\le n\le\langle h_i,\la\rangle\text{ and }u\in\text{Ker}\;e_i\cap M_{\la},\\
&\text{then }v=f_i^{(\langle h_i,\la\rangle)}u\text{ belongs to }\text{Ker}\;f_i\cap M_{s_i(\la)}
\text{ and }f_i^{(n)}u=e_i^{(\langle h_i,\la\rangle-n)}v.
\end{align*}
Here $s_i(\la)=\la-\langle h_i,\la\rangle\alpha_i$. Hence we obtain
\begin{equation} \label{Kop2}
\ft{i}(e_i^{(n)}v)=e_i^{(n-1)}v\quad\text{and}\quad \et{i}(e_i^{(n)}v)=e_i^{(n+1)}v
\end{equation}
for $v\in\text{Ker}\;f_i\cap M_{\mu}$ with $0\le n\le-\langle h_i,\mu\rangle$.

Let us now look at the definition of a crystal base. 
Let $A$ be the subring of $\Q(q)$ consisting of rational functions without poles at $q=0$.
Let $M$ be an integrable $\Uq$-module.

\begin{definition} \label{def:crystal base}
A pair $(L,B)$ is called a crystal base of $M$ if it satisfies the following 6 conditions:
\begin{eqnarray}
&&\text{$L$ is a free sub-$A$-module of $M$ such that $M\simeq\Q(q)\otimes_A L$,}\\
&&\text{$B$ is a base of the $\Q$-vector space $L/qL$,}\\
&&\et{i}L\subset L\text{ and }\ft{i}L\subset L\text{ for any }i. \label{eq:crystal base 3}
\end{eqnarray}
By (\ref{eq:crystal base 3}) $\et{i}$ and $\ft{i}$ act on $L/qL$.
\begin{eqnarray}
&&\et{i}B\subset B\cup\{0\}\text{ and }\ft{i}B\subset B\cup\{0\}.\\
&&L=\oplus_{\la\in P}L_{\la}\text{ and }B=\sqcup_{\la\in P}B_{\la}
\end{eqnarray}
where $L_{\la}=L\cap M_{\la}$ and $B_{\la}=B\cap(L_{\la}/qL_{\la})$.
\begin{equation} \label{eq:crystal base 6}
\text{For }b,b'\in B,b'=\ft{i}b\text{ if and only if }\et{i}b'=b.
\end{equation}
\end{definition}

Standard notations are in order. For $b\in B$ we set
\begin{align}
\veps_i(b)&=\max\{m\in\Zn\mid\et{i}^mb\neq0\}, &\vphi_i(b)&=\max\{m\in\Zn\mid\ft{i}^mb\neq0\},
\label{eq:epsphi1}\\
\veps(b)&=\sum_i\veps_i(b)\La_i, &\vphi(b)&=\sum_i\vphi_i(b)\La_i, \label{eq:epsphi2}\\
\wt b&=\vphi(b)-\veps(b). \label{eq:epsphi3}
\end{align}
Here $\{\La_i\}$ stands for the set of fundamental weights of $\geh$.

The crystal base behaves nicely under the tensor product. Let $(L_j,B_j)$ be the 
crystal base of an integrable $\Uq$-module $M_j$ ($j=1,2$). Set $L=L_1\ot_A L_2$
and $B=\{b_1\ot b_2\mid b_j\in B_j(j=1,2)\}$. Then $(L,B)$ is a crystal base of $M_1\ot M_2$.
Moreover, the action of $\et{i}$ and $\ft{i}$ becomes very simple as 
\begin{eqnarray}
\et{i}(b_1\ot b_2)&=&
\begin{cases}
\et{i}b_1\ot b_2&\text{ if }\vphi_i(b_1)\ge\veps_i(b_2),\\
b_1\ot \et{i}b_2&\text{ if }\vphi_i(b_1)<\veps_i(b_2),
\end{cases}\\
\ft{i}(b_1\ot b_2)&=&
\begin{cases}
\ft{i}b_1\ot b_2&\text{ if }\vphi_i(b_1)>\veps_i(b_2),\\
b_1\ot \ft{i}b_2&\text{ if }\vphi_i(b_1)\le\veps_i(b_2).
\end{cases}
\end{eqnarray}
Here $0\ot b$ and $b\ot0$ are understood to be $0$.
We denote this $B$ by $B_1\ot B_2$. $\veps_i,\vphi_i$ and $\wt$ are given by
\begin{eqnarray}
\veps_i(b_1\ot b_2)&=&
\max(\veps_i(b_1),\veps_i(b_1)+\veps_i(b_2)-\vphi_i(b_1)),\label{eq:ot-eps}\\
\vphi_i(b_1\ot b_2)&=&
\max(\vphi_i(b_2),\vphi_i(b_1)+\vphi_i(b_2)-\veps_i(b_2)),\label{eq:ot-phi}\\
\wt(b_1\ot b_2)&=&\wt b_1+\wt b_2.
\end{eqnarray}

Next lemma is used later in section \ref{sec:pr1}.

\begin{lemma} \label{lem:mod qL}
Let $(L,B)$ be a crystal base.
Assume that $\et{i}^3b=\ft{i}^3b=0$ for any $b\in B$. Let $v\in L$ 
be such that $v\equiv b$ mod $qL$. Then we have 
\begin{align*}
e_iv&\equiv q_i^{-\vphi_i(b)}\et{i}v\mbox{ mod }qq_i^{-\vphi_i(b)}L,\\
f_iv&\equiv q_i^{-\veps_i(b)}\ft{i}v\mbox{ mod }qq_i^{-\veps_i(b)}L.
\end{align*}
In particular, $e_iv\equiv0$ (resp. $f_iv\equiv0$) if $\veps_i(b)=0$ (resp. $\vphi_i(b)=0$).
\end{lemma}

\begin{proof}
We prove the second relation. Let $\la$ be the weight of $v$. From the assumption it suffices to
check the relation for the following cases, since the other cases are trivial.
\begin{itemize}
\item[(i)] $\veps_i(b)=0,\langle h_i,\la\rangle=1\text{ or }2$,
\item[(ii)] $\veps_i(b)=0\text{ or }1,\langle h_i,\la\rangle=0$.
\end{itemize}
In case (i) we have $f_iv=\ft{i}v$ by \eqref{sl2decomp1} and \eqref{Kop1}. In case (ii) let us
write $v=f_iv_1+v_2$ with $v_j\in \text{Ker}\;e_i\cap L$ ($j=1,2$) such that 
$\langle h_i,\wt v_1\rangle=2,\langle h_i,\wt v_2\rangle=0$. Then we have $f_iv=f_i^2v_1
=[2]_i\ft{i}v\equiv q_i^{-1}\ft{i}v$ mod $qq_i^{-1}L$.

The first relation can be checked similarly by using \eqref{sl2decomp2} and \eqref{Kop2}.
\end{proof}

\subsection{Polarization}

We define a total order on $\Q(q)$ by
\[
f>g\text{ if and only if }f-g\in\bigsqcup_{n\in\Z}\{q^n(c+qA)\mid c>0\}
\]
and $f\ge g$ if $f>g$ or $f=g$. 

Let $M$ and $N$ be $\Uq$-modules. A bilinear form $(\;,\;):M\ot_{\Q(q)}N\rightarrow\Q(q)$ is called 
an admissible pairing if it satisfies
\begin{align}
(q^hu,v)&=(u,q^hv),\nonumber\\
(e_iu,v)&=(u,q_i^{-1}t_i^{-1}f_iv), \label{admissible}\\
(f_iu,v)&=(u,q_i^{-1}t_ie_iv),\nonumber
\end{align}
for all $u\in M$ and $v\in N$. \eqref{admissible} implies
\begin{equation} \label{admissible2}
(e_i^{(n)}u,v)=(u,q_i^{-n^2}t_i^{-n}f_i^{(n)}v),\quad
(f_i^{(n)}u,v)=(u,q_i^{-n^2}t_i^{n}e_i^{(n)}v).
\end{equation}
A symmetric bilinear form $(\;,\;)$ on $M$ is called a 
preporlarization of $M$ if it satisfies \eqref{admissible} for $u,v\in M$.
A preporlarization is called a porlarization if it is positive definite with respective to 
the order on $\Q(q)$.

\subsection{Fusion construction} \label{subsec:fusion}

In what follows we assume that $\geh$ is of affine type. Let $P$ be the weight lattice, 
$\{\La_i\}$ the set of fundamental weights and $\delta$ the generator of null roots of $\geh$.
Then we have $P=\bigoplus_i\Z\La_i\oplus\Z\delta$. We set 
\[
P_{cl}=P/\Z\delta.
\]
Similar to the quantum algebra $\Uq$ which is associated with $P$, we can also consider 
$\Uqp$, which is associated with $P_{cl}$, namely, the subalgebra of $\Uq$ generated by 
$e_i,f_i,q^h$ ($h\in(P_{cl})^\ast$).

Let $K$ be a commutative ring containing $\Q(q)$ and let $x$ be an invertible element of $K$.
We introduce a $K\ot_{\Q(q)}\Uqp$-module $V_x$ by replacing the actions of $e_i,f_i$ with 
$x^{\delta_{i0}}e_i,x^{-\delta_{i0}}f_i$. The action of $q^h$ is not changed. Let $y$ also be 
an invertible element of $K$. A $K\ot_{\Q(q)}\Uqp$-linear map 
\[
R(x,y):\;V_x\ot V_y\longrightarrow V_y\ot V_x
\]
is called a $R$-matrix. Here we need to specify the coproduct $\Delta$ of $\Uq$ we use in this
paper. Our choice is the ``lower" one (see \cite{K1}) given by 
\begin{align*}
\Delta(q^h)&=q^h\ot q^h\quad\mbox{for }h\in(P_{cl})^\ast,\\
\Delta(e_i)&=e_i\ot t_i^{-1}+1\ot e_i,\\
\Delta(f_i)&=f_i\ot1+t_i\ot f_i.
\end{align*}
For a finite-dimensional $\Uqp$-module $V$ we assume the following.
\begin{align}
& \text{$V$ is irreducible}. \label{Virr}\\
& \text{There exists $\la_0\in P_{cl}$ such that 
 $\wt V\subset\la_0+\sum_{i\neq0}\Z_{\le0}\alpha_i\text{ and }\dim V_{\la_0}=1$}. \label{la0}
\end{align}
Here $\{\alpha_i\}$ is the set of simple roots. Under these assumptions it is known that there
exists a unique $R$-matrix up to a scalar multiple. Moreover, $R(x,y)$ depends only on $x/y$. 
Take a non zero vector $u_0$ from $V_{\la_0}$. We normalize $R(x,y)$ in such a way that
$R(x,y)(u_0\ot u_0)=u_0\ot u_0$. It is known in \cite{K3} that if $V$ is a ``good" module then 
the normalized $R$-matrix does not have a pole at $x/y=a\in A$. 

Next we review the fusion construction following section 3 of \cite{KMN2}.
Let $l$ be a positive integer and $\mathfrak{S}_l$ the $l$-th symmetric group.
Let $s_i$ be the simple reflection which interchanges $i$ and $i+1$, and 
let $\ell(w)$ be the length of $w\in\mathfrak{S}_l$. Let $V$ be
a finite-dimensional $U'_q(\geh)$-module. Let $R(x,y)$ denote the 
$R$-matrix for $V\ot V$. For any $w\in\mathfrak{S}_l$ we construct a map 
$R_w(x_1,\ldots,x_l):V_{x_1}\ot\cd\ot V_{x_l}\rightarrow V_{x_{w(1)}}\ot\cd\ot V_{x_{w(l)}}$
by
\begin{align*}
R_1(x_1,\ldots,x_l)&=1,\\
R_{s_i}(x_1,\ldots,x_l)&=\left(\bigotimes_{j<i}\text{id}_{V_{x_j}}\right)\ot
R(x_i,x_{i+1})\ot\left(\bigotimes_{j>i+1}\text{id}_{V_{x_j}}\right),\\
R_{ww'}(x_1,\ldots,x_l)&=R_{w'}(x_{w(1)},\ldots,x_{w(l)})\circ R_w(x_1,\ldots,x_l)\\
&\hspace{1cm}\text{ for $w,w'$ such that $\ell(ww')=\ell(w)+\ell(w')$.}
\end{align*}
Fix $r\in\Z_{>0}$. For each $l\in\Z_{>0}$, we put
\begin{align*}
R_l=&R_{w_0}(q^{r(l-1)},q^{r(l-3)},\ldots,q^{-r(l-1)}):\\
&V_{q^{r(l-1)}}\ot V_{q^{r(l-3)}}\ot\cd\ot V_{q^{-r(l-1)}}\rightarrow
V_{q^{-r(l-1)}}\ot V_{q^{-r(l-3)}}\ot\cd\ot V_{q^{r(l-1)}},
\end{align*}
where $w_0$ is the longest element of $\mathfrak{S}_l$. Then $R_l$ is a $U'_q(\geh)$-linear
homomorphism. Define 
\[
V_l=\mbox{Im}\;R_l.
\]
Let us denote by $W$ the image of 
\[
R(q^r,q^{-r}):V_{q^r}\ot V_{q^{-r}}\longrightarrow V_{q^{-r}}\ot V_{q^r}
\]
and by $N$ its kernel. Then we have
\begin{align}
&V_l\text{ considered as a submodule of }V^{\ot l}=V_{q^{-r(l-1)}}\ot\cd\ot V_{q^{r(l-1)}}\\
&\text{is contained in }\bigcap_{i=0}^{l-2}V^{\ot i}\ot W\ot V^{\ot(l-2-i)}. \nonumber\\
\intertext{Similarly, we have}
&V_l\text{ is a quotient of }V^{\ot l}/\sum_{i=0}^{l-2}V^{\ot i}\ot N\ot V^{\ot(l-2-i)}.
\label{Vlquotient}
\end{align}

\subsection{Preliminary propositions}

In this subsection, following \cite{KMN2} we define a prepolarization on $V_l$ and prepare
a necessary proposition to show the main theorem. First we recall

\begin{lemma}
Let $M_j$ and $N_j$ be $\Uqp$-modules and let $(\;,\;)_j$ be an admissible pairing between 
$M_j$ and $N_j$ $(j=1,2)$. Then the pairing $(\;,\;)$between $M_1\ot M_2$ and $N_1\ot N_2$
defined by $(u_1\ot u_2,v_1\ot v_2)=(u_1,v_1)_1(u_2,v_2)_2$ for all $u_j\in M_j$ and 
$v_j\in N_j$ is admissible.
\end{lemma}

Let $V$ be a finite-dimensional $\Uqp$-module satisfying \eqref{Virr},\eqref{la0}. Suppose
$V$ has a polarization.
The polarization on $V$ gives an admissible pairing between $V_x$ and $V_{x^{-1}}$. Hence
it induces an admissible pairing between $V_{x_1}\ot\cd\ot V_{x_l}$ and 
$V_{x_1^{-1}}\ot\cd\ot V_{x_l^{-1}}$.

\begin{lemma}
If $x_j=x_{l+1-j}^{-1}$ for $j=1,\ldots,l$, then for any $u,u'\in V_{x_1}\ot\cd\ot V_{x_l}$,
we have
\[
(u,R_{w_0}(x_1,\ldots,x_l)u')=(u',R_{w_0}(x_1,\ldots,x_l)u).
\]
\end{lemma}
By taking $x_1=q^{r(l-1)},x_2=q^{r(l-3)}$, etc., we obtain the admissible pairing $(\;,\;)$
between $W=V_{q^{r(l-1)}}\ot V_{q^{r(l-3)}}\ot\cd\ot V_{q^{-r(l-1)}}$ and 
$W'=V_{q^{-r(l-1)}}\ot V_{q^{-r(l-3)}}\ot\cd\ot V_{q^{r(l-1)}}$ that satisfies
\begin{equation} \label{wRw'}
(w,R_lw')=(w',R_lw)\quad\text{for any }w,w'\in W.
\end{equation}
This allows us to define a preporlarization $(\;,\;)_l$ on $V_l$ by
\[
(R_lu,R_lu')_l=(u,R_lu')
\]
for $u,u'\in V_{q^{r(l-1)}}\ot V_{q^{r(l-3)}}\ot\cd\ot V_{q^{-r(l-1)}}$. 

Next we introduce a $\Z$-form of $\Uqp$. 
Recall that $A$ is the subring of $\Q(q)$ consisting of 
rational functions without poles at $q=0$. We introduce the subalgebras $A_{\Z}$ and $\Kz$ 
of $\Q(q)$ by
\begin{align*}
A_{\Z}&=\{f(q)/g(q)\mid f(q),g(q)\in\Z[q],g(0)=1\},\\
\Kz&=A_{\Z}[q^{-1}].
\end{align*}
Then we have 
\[
\Kz\cap A=A_{\Z},\quad A_{\Z}/qA_{\Z}\simeq\Z.
\]
We then define $\Uqp_{K_{\Z}}$ as the $K_{\Z}$-subalgebra of $\Uqp$ generated by 
$e_i,f_i,q^h$ ($h\in(P_{cl})^\ast$). Set $V_{K_{\Z}}=\Uqp_{K_{\Z}}u_0$ and assume
\begin{equation} \label{Kzu0}
(V_{K_{\Z}})_{\la_0}=K_{\Z}u_0.
\end{equation}
Let us further set
\[
(V_l)_{K_{\Z}}=R_l((V_{K_{\Z}})^{\ot l})\cap(V_{K_{\Z}})^{\ot l}.
\]
Then one can show

\begin{proposition}
\begin{itemize}
\item[(i)] $(\;,\;)_l$ is a nondegenerate prepolarization on $V_l$.
\item[(ii)] $(R_l(u_0^{\ot l}),R_l(u_0^{\ot l}))_l=1$.
\item[(iii)] $((V_l)_{\Kz},(V_l)_{\Kz})_l\subset\Kz$.
\end{itemize}
\end{proposition}

Let $I$ be the index set of the vertices of the Dynkin diagram of $\geh$ with the vertex $0$
as in \cite{Kac}. Let $\geh_0$ be the finite-dimensional simple Lie algebra whose Dynkin 
diagram is obtained by removing the $0$-vertex from that of $\geh$. Let $\ol{P}_+$ be the set
of dominant integral weights of $\geh_0$ and $V(\la)$ be the irreducible highest weight 
$U_q(\geh_0)$-module of highest weight $\la$ for $\la\in\ol{P}_+$. The following proposition,
which is essentially stated in Proposition 2.6.1 and 2.6.2 of \cite{KMN2},
is a key to prove the main theorem.

\begin{proposition} \label{prop:key}
Let $M$ be a finite-dimensional integrable $\Uqp$-module. Let $(\;,\;)$ be a prepolarization
on $M$, and $M_{\Kz}$ a $\Uqp_{K_{\Z}}$-submodule of $M$ such that $(M_{\Kz},M_{\Kz})\subset\Kz$. 
Let $\la_1,\ldots,\la_m\in\ol{P}_+$, and we assume the following conditions.
\begin{equation} \label{char}
\dim M_{\la_k}\le\sum_{j=1}^m\dim V(\la_j)_{\la_k}\mbox{ for }k=1,\ldots,m.
\end{equation}
\begin{align}
&\text{There exist $u_j\in(M_{\Kz})_{\la_j}$ $(j=1,\ldots,m)$ such that 
$(u_j,u_k)\in\delta_{jk}+qA$,} \label{norm}\\
&\text{and $(e_iu_j,e_iu_j)\in qq_i^{-2(1+\langle h_i,\la_j\rangle)}A$ for any $i\in I\setminus\{0\}$.}
\nonumber
\end{align}
Set $L=\{u\in M\mid(u,u)\in A\}$ and set $B=\{b\in M_{\Kz}\cap L/M_{\Kz}\cap qL\mid(b,b)_0=1\}$.
Here $(\;,\;)_0$ is the $\Q$-valued symmetric bilinear form on $L/qL$ induced by $(\;,\;)$.
Then we have the following.
\begin{itemize}
\item[(i)] $(\;,\;)$ is a polarization on $M$.
\item[(ii)] $M\simeq\bigoplus_j V(\la_j)$ as $U_q(\geh_0)$-modules.
\item[(iii)] $(L,B)$ is a crystal pseudobase of $M$.
\end{itemize}
\end{proposition}

\section{KR module of type $D$}

\subsection{KR module $W_1^{(k)}$}

First we review the Dynkin datum of type $D_n^{(1)}$. Let $I=\{0,1,\ldots,n\}$ be the index set
of the Dynkin diagram, $\{\alpha_i\}_{i\in I}$ the set of simple roots, $\{\La_i\}_{i\in I}$ the 
set of fundamental weights. The standard null root $\delta$ is given by 
\begin{equation} \label{null root}
\delta=\alpha_0+\alpha_1+2\alpha_2+\cd+2\alpha_{n-2}+\alpha_{n-1}+\alpha_n.
\end{equation}
We denote the weight lattice by $P$, that is, $P=\bigoplus_{i\in I}\Z\La_i\oplus\Z\delta$.
The sublattice $\ol{P}=\bigoplus_{i\in I_0}\Z\olLa_i$ can be viewed as the weight lattice for 
$D_n$. Here $I_0=I\setminus\{0\}$ and $\olLa_i=\La_i-a_i\La_0$ with $a_i$ being the coefficient
of $\alpha_i$ in \eqref{null root}. It is sometimes useful to introduce an orthonormal basis 
$\{\eps_1,\eps_2,\ldots,\eps_n\}$ of $\Q\ot_{\Z}\ol{P}$ in such a way that we have
\begin{align*}
\alpha_i&=\left\{
\begin{array}{ll}
\eps_i-\eps_{i+1}\quad&(i=1,\ldots,n-1)\\
\eps_{n-1}+\eps_n&(i=n),
\end{array}\right.\\
\olLa_i&=\left\{
\begin{array}{ll}
\eps_1+\cd+\eps_i&(i=1,\ldots,n-2)\\
(\eps_1+\cd+\eps_{n-1}-\eps_n)/2\quad&(i=n-1)\\
(\eps_1+\cd+\eps_{n-1}+\eps_n)/2\quad&(i=n).
\end{array}\right.
\end{align*}
Then we have $\alpha_0=\delta-\eps_1-\eps_2$.
Since the lengths of the simple roots are all equal, we have $q_i=q$ for any $i\in I$.
Hence we shall abbreviate $i$ from $[m]_i$ or $[m]_i!$.

Let $W_1^{(k)}$ be the $k$-th fundamental representation of $U'_q(D_n^{(1)})$. It is known that
it has the following decomposition into $U_q(D_n)$-modules.
\begin{equation} \label{decomp}
W_1^{(k)}\simeq\left\{
\begin{array}{ll}
V(\olLa_k)\oplus V(\olLa_{k-2})\oplus\cd\oplus V(\olLa_1\text{ or }0)\;&\text{if }1\le k\le n-2,\\
V(\olLa_k)&\text{if }k=n-1,n.
\end{array}\right.
\end{equation}
On $W_1^{(k)}$ the following results are known.
\begin{proposition} \label{prop:W_1^(k)}
\begin{itemize}
\item[(1)] $W_1^{(k)}$ is ``good" in the sense of Kashiwara. In particular, it has a crystal base.
\item[(2)] $W_1^{(k)}$ has a polarization.
\end{itemize}
\end{proposition}
The first claim is due to Kashiwara \cite{K3} and the second to Koga \cite{Ko}, who got the result 
by exploiting the fusion construction among the spin representations.

\subsection{Crystal of $W_1^{(k)}$}

We denote the crystal of $W_1^{(k)}$ by $B^{k,1}$. We review in this subsection Schilling's 
variation of Koga's result on the crystal structure of $B^{k,1}$. First we treat the case of
$k=1$. The crystal graph of $B^{1,1}$ is depicted as follows.

\vspace{5mm}\noindent
\unitlength 0.1in
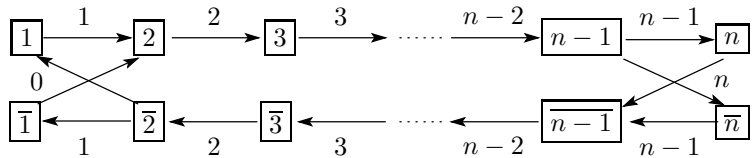
\begin{picture}( 43.2800,  6.8600)( -2.5100,-12.3900)
%
\special{pn 8}%
\special{pa 1194 750}%
\special{pa 1646 750}%
\special{fp}%
\special{sh 1}%
\special{pa 1646 750}%
\special{pa 1580 730}%
\special{pa 1594 750}%
\special{pa 1580 770}%
\special{pa 1646 750}%
\special{fp}%
%
\special{pn 8}%
\special{pa 522 750}%
\special{pa 976 750}%
\special{fp}%
\special{sh 1}%
\special{pa 976 750}%
\special{pa 908 730}%
\special{pa 922 750}%
\special{pa 908 770}%
\special{pa 976 750}%
\special{fp}%
%
\special{pn 8}%
\special{pa 1858 750}%
\special{pa 2310 750}%
\special{fp}%
\special{sh 1}%
\special{pa 2310 750}%
\special{pa 2242 730}%
\special{pa 2256 750}%
\special{pa 2242 770}%
\special{pa 2310 750}%
\special{fp}%
%
\special{pn 8}%
\special{pa 2656 750}%
\special{pa 3110 750}%
\special{fp}%
\special{sh 1}%
\special{pa 3110 750}%
\special{pa 3042 730}%
\special{pa 3056 750}%
\special{pa 3042 770}%
\special{pa 3110 750}%
\special{fp}%
%
\special{pn 8}%
\special{pa 3570 760}%
\special{pa 4024 760}%
\special{fp}%
\special{sh 1}%
\special{pa 4024 760}%
\special{pa 3956 740}%
\special{pa 3970 760}%
\special{pa 3956 780}%
\special{pa 4024 760}%
\special{fp}%
%
\special{pn 8}%
\special{pa 3110 1190}%
\special{pa 2656 1190}%
\special{fp}%
\special{sh 1}%
\special{pa 2656 1190}%
\special{pa 2724 1210}%
\special{pa 2710 1190}%
\special{pa 2724 1170}%
\special{pa 2656 1190}%
\special{fp}%
%
\special{pn 8}%
\special{pa 4050 1190}%
\special{pa 3598 1190}%
\special{fp}%
\special{sh 1}%
\special{pa 3598 1190}%
\special{pa 3664 1210}%
\special{pa 3650 1190}%
\special{pa 3664 1170}%
\special{pa 3598 1190}%
\special{fp}%
%
\special{pn 8}%
\special{pa 2310 1190}%
\special{pa 1858 1190}%
\special{fp}%
\special{sh 1}%
\special{pa 1858 1190}%
\special{pa 1924 1210}%
\special{pa 1910 1190}%
\special{pa 1924 1170}%
\special{pa 1858 1190}%
\special{fp}%
%
\special{pn 8}%
\special{pa 1638 1190}%
\special{pa 1186 1190}%
\special{fp}%
\special{sh 1}%
\special{pa 1186 1190}%
\special{pa 1254 1210}%
\special{pa 1240 1190}%
\special{pa 1254 1170}%
\special{pa 1186 1190}%
\special{fp}%
%
\special{pn 8}%
\special{pa 976 1186}%
\special{pa 522 1186}%
\special{fp}%
\special{sh 1}%
\special{pa 522 1186}%
\special{pa 590 1206}%
\special{pa 576 1186}%
\special{pa 590 1166}%
\special{pa 522 1186}%
\special{fp}%
\put(4.3200,-7.5000){\makebox(0,0){\fbox{$1$}}}%
\put(10.7300,-7.5000){\makebox(0,0){\fbox{$2$}}}%
%
\special{pn 8}%
\special{pa 3558 862}%
\special{pa 4078 1100}%
\special{fp}%
\special{sh 1}%
\special{pa 4078 1100}%
\special{pa 4026 1054}%
\special{pa 4030 1078}%
\special{pa 4008 1090}%
\special{pa 4078 1100}%
\special{fp}%
%
\special{pn 8}%
\special{pa 4078 862}%
\special{pa 3558 1108}%
\special{fp}%
\special{sh 1}%
\special{pa 3558 1108}%
\special{pa 3626 1098}%
\special{pa 3606 1084}%
\special{pa 3610 1060}%
\special{pa 3558 1108}%
\special{fp}%
%
\special{pn 8}%
\special{pa 492 1096}%
\special{pa 1014 858}%
\special{fp}%
\special{sh 1}%
\special{pa 1014 858}%
\special{pa 944 868}%
\special{pa 966 880}%
\special{pa 962 904}%
\special{pa 1014 858}%
\special{fp}%
%
\special{pn 8}%
\special{pa 1014 1096}%
\special{pa 492 850}%
\special{fp}%
\special{sh 1}%
\special{pa 492 850}%
\special{pa 544 898}%
\special{pa 540 874}%
\special{pa 562 860}%
\special{pa 492 850}%
\special{fp}%
\put(17.5900,-7.5000){\makebox(0,0){\fbox{$3$}}}%
\put(33.3500,-7.4700){\makebox(0,0){\fbox{$n-1$}}}%
\put(41.2500,-7.5500){\makebox(0,0){\fbox{$n$}}}%
%
\special{pn 8}%
\special{pa 2386 750}%
\special{pa 2604 750}%
\special{dt 0.045}%
%
\special{pn 8}%
\special{pa 2386 1190}%
\special{pa 2604 1190}%
\special{dt 0.045}%
\put(4.2400,-11.9000){\makebox(0,0){\fbox{$\ol{1}$}}}%
\put(10.7300,-11.9000){\makebox(0,0){\fbox{$\ol{2}$}}}%
\put(17.3600,-11.9000){\makebox(0,0){\fbox{$\ol{3}$}}}%
\put(33.3500,-11.8700){\makebox(0,0){\fbox{$\ol{n-1}$}}}%
\put(41.2500,-11.9500){\makebox(0,0){\fbox{$\ol{n}$}}}%
\put(7.3300,-6.3800){\makebox(0,0){$1$}}%
\put(14.1200,-6.3800){\makebox(0,0){$2$}}%
\put(20.7600,-6.3800){\makebox(0,0){$3$}}%
\put(28.7500,-6.3800){\makebox(0,0){$n-2$}}%
\put(37.9800,-6.4500){\makebox(0,0){$n-1$}}%
\put(4.8500,-9.7700){\makebox(0,0){$0$}}%
\put(40.7000,-9.8000){\makebox(0,0){$n$}}%
\put(7.3300,-13.1700){\makebox(0,0){$1$}}%
\put(14.1200,-13.1700){\makebox(0,0){$2$}}%
\put(20.7600,-13.1700){\makebox(0,0){$3$}}%
\put(28.7500,-13.1700){\makebox(0,0){$n-2$}}%
\put(37.9800,-13.2400){\makebox(0,0){$n-1$}}%
\end{picture}%

\vspace{8mm}\noindent
Here for $b,b'\in B^{1,1}$ $b\stackrel{i}{\longrightarrow}$ means $\ft{i}b=b'$ ($\Leftrightarrow
b=\et{i}b'$).

Next in view of \eqref{decomp} we recall the $U_q(D_n)$-crystal structure of $B(\olLa_l)$, the
crystal of $U_q(D_n)$-module $V(\olLa_l)$, by \cite{KN}. Consider the alphabet $\mathcal{A}=
\{1,2,\ldots,n,\ol{n},\ol{n-1},$ $\ldots,\ol{1}\}$ consisting of the crystal elements of $B^{1,1}$.
It is given the following (partial) order.
\[
1\prec2\prec\cd\prec n-1\prec{n\atop\ol{n}}\prec\ol{n-1}\prec\cd\prec\ol{1}.
\]
Then, for $1\le l\le n-2$ $B(\olLa_l)$ is identified with the set of columns 
\begin{center}
\begin{tabular}{|c|}
\hline $m_1$\\ \hline $m_2$\\ \hline $\vdots$\\ \hline $m_l$\\ \hline
\end{tabular}
\end{center}
of height $l$ satisfying 
\begin{align}
&m_j\not\succ m_{j+1}\quad\text{for }j=1,\ldots,l-1, \label{tab cond1}\\
&\text{if $m_a=p$ and $m_b=\ol{p}$, then dist$(p,\ol{p})\le p$.} \label{tab cond2}
\end{align}
Here dist$(p,\ol{p})=a+l+1-b$ if $m_a=p$ and $m_b=\ol{p}$. The column tableau as above is also
written as $m_1m_2\cd m_l$. Note that we allow $(m_j,m_{j+1})=(n,\ol{n})$ and $(\ol{n},n)$.
The actions of $\et{i},\ft{i}$ ($i\in I_0$) is given by considering the embedding
\[
\begin{array}{ccc}
B(\olLa_l)&\hookrightarrow&(B^{1,1})^{\ot l}\\
m_1m_2\cd m_l&\mapsto&m_1\ot m_2\ot\cd\ot m_l
\end{array}
\]
and apply the tensor product rule of crystals on the r.h.s. $B(0)$ is realized as $\{\phi\}$ with 
the trivial actions of $\et{i},\ft{i}$ ($i\in I_0$), that is, $\et{i}\phi=\ft{i}\phi=0$.

For $1\le k\le n-2$, we are to represent $B^{k,1}$ as the set of column tableaux of height 
$k$ satisfying \eqref{tab cond1}. By \eqref{decomp} $B^{k,1}$ is the union of the sets 
corresponding to $B(\olLa_l)$ with $0\le l\le k$ and $l\equiv k$ (mod 2). In \cite{S} maps 
from $B(\olLa_l)$ to column tableaux of height $k$ were defined. If $b\in
B(\olLa_l)$, then fill the column of height $l$ of $b$ succesively by a pair 
$(i_j,\ol{i_j})$ for $1\le j\le(k-l)/2$ in the following way to obtain a column of height $k$.
Set $i_0=0$. Let $i_{j-1}<i_j$ be minimal such that
\begin{itemize}
\item[(1)] neither $i_j$ or $\ol{i_j}$ is in the column;
\item[(2)] adding $i_j$ and $\ol{i_j}$ to the column we have dist$(i_j,\ol{i_j})\ge i_j+j$;
\item[(3)] adding $i_j$ and $\ol{i_j}$ to the column, all other pairs $(a,\ol{a})$ in the new 
	column with $a>i_j$ satisfy dist$(a,\ol{a})\le a+j$.
\end{itemize}
The filling map and $\ft{i}$ for $i\in I_0$ commute. Denote the filling map to height $k$ by
$F_k$ or simply $F$. Let $D_k$ or $D$, the dropping map, be the inverse of $F_k$. 
Explicitly, given a one-column tableau
of $b$ of height $k$, let $i_0=0$ and successively find $i_j>i_{j-1}$ minimal such that the pair 
$(i_j,\ol{i_j})$ is in $b$ and dist$(i_j,\ol{i_j})\ge i_j+j$. Drop all such pairs $(i_j,\ol{i_j})$
from $b$. Thus we have 
\[
B^{k,1}\simeq\bigoplus_{0\le l\le k,\:l\equiv k\,(2)}F_k(B(\olLa_l))
\quad\mbox{as $U_q(D_n)$-crystals}.
\]
It is the set of all column tableaux of height $k$ satisfying \eqref{tab cond1} only.

We are left to give the rule of the actions of $\et{0}$ and $\ft{0}$.
For this purpose we need slight variants of $F_k$ and $D_k$, denoted by $\tilde{F}_k$
and $\tilde{D}_k$, respectively, which act on columns that do not contain $1,2,\ol{2},\ol{1}$.
On these columns $\tilde{F}_k$ and $\tilde{D}_k$ are defined by replacing $i\mapsto i-2$ and 
$\ol{i}\mapsto\ol{i-2}$, then applying $F_k$ and $D_k$, and finally replacing $i\mapsto i+2$ and 
$\ol{i}\mapsto\ol{i+2}$. The following proposition is given in \cite{S}.

\begin{proposition} \label{prop:e0f0}
For $b\in B^{k,1}$,
\[
\et{0}b=\left\{
\begin{array}{ll}
F_k(\tilde{D}_{k-2}(x))&\text{if }b=12x\\
\tilde{F}_{k-1}(x)\ol{2}&\text{if }b=12x\ol{2}\\
\tilde{F}_{k-1}(x)\ol{1}&\text{if }b=12x\ol{1}\\
\tilde{F}_{k-2}(x)\ol{2}\ol{1}&\text{if }b=12x\ol{2}\ol{1}\\
F_k(\tilde{D}_{k-1}(x)\ol{2})&\text{if }b=1x\\
F_k(\tilde{D}_{k-1}(x)\ol{1})&\text{if }b=2x\\
x\ol{2}\ol{1}&\text{if }b=1x\ol{1}\text{ and }\tilde{D}_{k-2}(x)=x\\
0&\text{otherwise}
\end{array}\right.
\]
\[
\ft{0}b=\left\{
\begin{array}{ll}
F_k(\tilde{D}_{k-2}(x))&\text{if }b=x\ol{2}\ol{1}\\
2\tilde{F}_{k-1}(x)&\text{if }b=2x\ol{2}\ol{1}\\
1\tilde{F}_{k-1}(x)&\text{if }b=1x\ol{2}\ol{1}\\
12\tilde{F}_{k-2}(x)&\text{if }b=12x\ol{2}\ol{1}\\
F_k(2\tilde{D}_{k-1}(x))&\text{if }b=x\ol{1}\\
F_k(1\tilde{D}_{k-1}(x))&\text{if }b=x\ol{2}\\
12x&\text{if }b=1x\ol{1}\text{ and }\tilde{D}_{k-2}(x)=x\\
0&\text{otherwise}
\end{array}\right.
\]
where $x$ does not contain $1,2,\ol{2},\ol{1}$.
\end{proposition}

\subsection{Existence of crystal pseudobase for $W_l^{(k)}$}

In this subsection we prove our main theorem by using Proposition \ref{prop:key}. We prepare
several lemmas and propositions.

\begin{lemma}
Let $R(x/y)$ be the $R$-matrix from $W_{1,x}^{(k)}\ot W_{1,y}^{(k)}$ to 
$W_{1,y}^{(k)}\ot W_{1,x}^{(k)}$. Then it has the following form.
\[
R(z)=P_{2\varpi_k}+\frac{z-q^2}{1-q^2z}P_{\varpi_{k+1}+\varpi_{k-1}}+\cd.
\]
Here $z=x/y,\varpi_j=\eps_1+\eps_2+\cd+\eps_j\in\ol{P}_+$ for $0\le j\le n-1$, and $P_{\la}$
stands for the projector onto the irreducible $U_q(D_n)$-module $V(\la)$ in $(W_1^{(k)})^{\ot2}$.
\end{lemma}

\begin{proof}
Let $u_0$ be a $U_q(D_n)$-highest weight vector of $W_1^{(k)}$ of highest weight 
$\olLa_k(=\varpi_k)$. Since 
$V(\varpi_{k+1}+\varpi_{k-1})$ is multiplicity free, a unique highest weight vector up to a 
scalar is given by
\[
v=u_0\ot f_ku_0-qf_ku_0\ot u_0.
\]
Since $f_iu_0=0$ for $i\in I\setminus\{k\}$, we have
\[
F^{(1)}v=u_0\ot F^{(1)}u_0-qF^{(1)}u_0\ot u_0
\]
where $F^{(1)}=f_{k+1}\cd f_{n-2}f_nf_{n-1}\cd f_{k+1}f_1\cd f_{k-1}$. Hence we have
\[
F^{(2)}v=qu_0\ot F^{(2)}u_0-qF^{(2)}u_0\ot u_0+(\text{unwanted terms})
\]
where $F^{(2)}=f_2\cd f_kF^{(1)}$ and we know $f_0(\text{unwanted terms})=0$ by weight consideration.
Hence we have
\[
F^{(3)}v=q^{-1}y^{-1}u_0\ot F^{(3)}u_0-qx^{-1}F^{(3)}u_0\ot u_0\text{ on }V_x\ot V_y
\]
where $F^{(3)}=f_0F^{(2)}$ and $V=W_1^{(k)}$. Note that $F^{(3)}u_0=\alpha u_0$ with some 
$\alpha\neq0$, since the corresponding crystal element is not killed from Proposition 
\ref{prop:e0f0}. Thus we have
\[
F^{(3)}v=\alpha(q^{-1}y^{-1}-qx^{-1})u_0\ot u_0.
\]
Now let
\[
R(z)\propto\vphi(z)P_{2\varpi_k}+\vphi'(z)P_{\varpi_{k+1}+\varpi_{k-1}}+\cd.
\]
Then we have
\begin{align*}
&R(z)F^{(3)}v=\alpha(q^{-1}y^{-1}-qx^{-1})\vphi(z)(u_0\ot u_0)\\
=&F^{(3)}R(z)v=\vphi'(z)F^{(3)}v'=\alpha(q^{-1}x^{-1}-qy^{-1})\vphi'(z)(u_0\ot u_0).
\end{align*}
Here by $v'$ we mean that it is considered to be in $V_y\ot V_x$. Thus we have
\[
\vphi'(z)/\vphi(z)=\frac{z-q^2}{1-q^2z}.
\]
\end{proof}

Set $W=\mbox{Im}\;R(q^2),N=\mbox{Ker}\;R(q^2)$. They are $U'_q(D_n^{(1)})$-modules. Using the main
result of \cite{N} one can show the following.

\begin{lemma} \label{lem:WN}
We have 
\[
W\simeq W^{(k)}_2,\quad N\simeq\bigotimes_{j\sim k}W_1^{(j)}
\]
as $U'_q(D_n^{(1)})$-modules. Here $j\sim k$ means that the corresponding vertices are tied by
an edge in the Dynkin diagram.
Moreover, both $W$ and $N$ are irreducible.
\end{lemma}

\begin{proof}
In \cite{N} it is shown that there exists an exact sequence of $U'_q(D_n^{(1)})$-modules
\[
0\longrightarrow\bigotimes_{j\sim k}W_1^{(j)}\longrightarrow W_{1,q}^{(k)}\ot W_{1,q^{-1}}^{(k)}
\longrightarrow W_2^{(k)}\longrightarrow0.
\]
(An acute reader should have noticed that the exact sequence is different from \cite{N}. It is 
because the definition of the KR modules and the choice of the coproduct are different.)
Moreover, it is also known that $\bigotimes_{j\sim k}W_1^{(j)}$ and $W_2^{(k)}$ are irreducible.
Set $W'=\bigotimes_{j\sim k}W_1^{(j)}$ and consider $N\cap W'$. Since $W'$ is irreducible, we have
$N\cap W'=\{0\}$ or $W'$. Recall that $W_{1,q}^{(k)}\ot W_{1,q^{-1}}^{(k)}$ contains a unique 
irreducible $U_q(D_n)$-module $V(\varpi_{k+1}+\varpi_{k-1})$. From the previous lemma and 
\eqref{decomp} we know it is contained both in $N$ and in $W'$. Hence we have $N\supset W'$. 
Now suppose $N\supsetneqq W'$. Then we have a surjective $U'_q(D_n^{(1)})$-linear map
\[
W_{1,q}^{(k)}\ot W_{1,q^{-1}}^{(k)}/W'\longrightarrow W_{1,q}^{(k)}\ot W_{1,q^{-1}}^{(k)}/N.
\]
Since the l.h.s. is irreducible, $N=W'\mbox{ or }W_{1,q}^{(k)}\ot W_{1,q^{-1}}^{(k)}$.
Since $N$ cannot be the second choice by the previous lemma. One obtains $N=W'$ and 
$W\simeq W_2^{(k)}$.
\end{proof}

Since $W$ is known to be a KR module by the previous lemma, we have
\begin{lemma} \label{lem:decomp of W}
As a $U_q(D_n)$-module $W$ has the following decomposition.
\[
W\simeq\bigoplus_{0\le m_1\le m_2\le [k/2]}V(\olLa_{k-2m_1}+\olLa_{k-2m_2})
\] 
\end{lemma}

We set $B=B^{k,1}$. We fix a basis $\{v_I\}_{I\in B}$ of $W_1^{(k)}$ in such a way that 
$v_I\mbox{ mod }qL=I$ as an element of $B$.

\begin{proposition} \label{prop:N}
$N$ contains a vector of the form
\[
v_{I_1}\ot v_{I_2}-\sum_{J_1\ot J_2\in B_1}a_{J_1J_2}v_{J_1}\ot v_{J_2}\quad(a_{J_1J_2}\in A)
\]
for any $I_1,I_2$ such that $I_1\ot I_2\in B^{\ot2}\setminus B_1$.
\end{proposition}
\noindent
See \eqref{B1h},\eqref{Ba} for the definition of $B_1$.

We now apply the fusion construction in section \ref{subsec:fusion} to $V=W_1^{(k)}$ with $r=1$.
The assumptions \eqref{Virr},\eqref{la0} are valid with $\la_0=\olLa_k$. \eqref{Kzu0} can also be 
checked. Other necessary properties are guaranteed by Proposition \ref{prop:W_1^(k)}.
For $l\in\Z_{>0}$ we define $W_l^{(k)}=\mbox{Im}\;R_l$. Let $k'=[k/2]$. 
Let $\cb=(c_1,c_2,\ldots,c_{k'})$ be a sequence of integers such that 
$l\ge c_1\ge c_2\ge\cd\ge c_{k'}\ge0$. For such $\cb$ we define a vector $u_m$ ($0\le m\le k'$) 
in $W_l^{(k)}$ inductively by
\begin{equation} \label{u_m}
u_m=(e_{k-2m}^{(c_m)}\cd e_2^{(c_m)}e_1^{(c_m)})(e_{k-2m+1}^{(c_m)}\cd e_3^{(c_m)}e_2^{(c_m)})
e_0^{(c_m)}u_{m-1},
\end{equation}
where $u_0$ here is $u_0^{\ot l}$ in $(W_1^{(k)})^{\ot l}$. Set $u(\cb)=u_{k'}$. The weight of 
$u(\cb)$ is given by
\[
\la(\cb)=\sum_{j=0}^{k'}(c_j-c_{j+1})\olLa_{k-2j},
\]
where we have set $c_0=l,c_{k'+1}=0$. 
For $l,m\in\Z_{\ge0}$ such that $m\le l$ we define the $q$-binomial coefficient by
\[
{l\brack m}=\frac{[l]!}{[m]![l-m]!}.
\]
The following proposition calculates values of prepolarizations 
necessary to prove the main theorem.
\begin{proposition} \label{prop:norm}
\begin{itemize}
\item[(1)] ${\displaystyle (u(\cb),u(\cb))_l=\prod_{j=1}^{k'}q^{c_j(2l-c_j)}{2l\brack c_j}}$,
\item[(2)] $(e_ju(\cb),e_ju(\cb))_l=0$ unless $k-j\in 2\Z_{\ge0}$. If $k-j\in 2\Z_{\ge0}$, then 
	setting $p=(k-j)/2+1$ $(e_ju(\cb),e_ju(\cb))_l$ is given by
\[
q^{2l-c_{p-1}-1}[2l-c_{p-1}]\prod_{j=1}^{k'}q^{(c_j-\delta_{j,p})(2l-c_j)}
{2l-\delta_{j,p}\brack c_j-\delta_{j,p}}.
\]
\end{itemize}
\end{proposition}
\noindent
Proofs of these propositions are given in subsequent sections.

The rest of this section is devoted to the proof of Theorem \ref{th:main}.
{}From Proposition \ref{prop:key} it suffices to show
\begin{itemize}
\item[(i)] $\mbox{\sl ch}\:W_l^{(k)}\le\sum_{l\ge c_1\ge\cd\ge c_{k'}\ge0}
\mbox{\sl ch}\:V(\la(\cb))$, where $V(\la)$ is the irreducible $U_q(D_n)$-module
	with highest weight $\la$ and $\mbox{\sl ch}\:V$ stands for the formal character of $V$.
\item[(ii)] $(u(\cb),u(\cb'))_l\in\delta_{\cb\cb'}+qA$ and 
	$(e_ju(\cb),e_ju(\cb))_l\in q^{-1-2\langle h_j,\la(\cb)\rangle}A$ for $j\ne0$.
\end{itemize}

Let us show (i). First notice that $\sum_{l\ge c_1\ge\cd\ge c_{k'}\ge0}\mbox{\sl ch}\:V(\la(\cb))
=\sum_{0\le m_1\le\cd\le m_l\le k'}$ $\mbox{\sl ch}\:V(\olLa_{k-2m_1}+\cd+\olLa_{k-2m_l})$.
In view of \eqref{Vlquotient} and Proposition \ref{prop:N} $W_l^{(k)}$ is a 
quotient of a module generated by the set of vectors
\[
\{v_{I_1}\ot v_{I_2}\ot\cd\ot v_{I_l}\mid I_j\ot I_{j+1}\in B_1\text{ for }j=1,\ldots,l-1\}.
\]
Assume $I_j\in B(\olLa_{k-2m_{l+1-j}})\subset B$. Then $I_1\ot I_2\ot\cd\ot I_l$ belongs to 
\begin{align*}
\bigcap_{j=1}^{l-1}B(\olLa_{k-2m_l})\ot\cd\ot B(\olLa_{k-2m_{j+2}})
&\ot B(\olLa_{k-2m_{j+1}}+\olLa_{k-2m_j})\\
&\ot B(\olLa_{k-2m_{j-1}})\ot\cd\ot B(\olLa_{k-2m_1}).
\end{align*}
However, the above crystal is known to be identified with $B(\olLa_{k-2m_1}+\cd+\olLa_{k-2m_l})$
(see Proposition 2.2.1 of \cite{KN}). This fact verifies (i).

For the proof of (ii) note that 
\[
[m]\in q^{1-m}A,\quad{m\brack n}\in q^{-n(m-n)}A.
\]
If $\cb\neq\cb'$, $(u(\cb),u(\cb'))_l=0$ since the weights of $u(\cb)$ and $u(\cb')$ are different.
$(u(\cb),u(\cb'))_l\in1+qA$ by Proposition \ref{prop:norm} (1). For the second part it suffices to 
notice that $\langle h_i,\la(\cb)\rangle\ge0$.
The proof is completed.

\section{Proof of proposition \ref{prop:N}} \label{sec:pr1}

We prepare several lemmas. The next one is a direct consequence of Proposition \ref{prop:e0f0}.

\begin{lemma} \label{lem:e0}
Suppose $k\ge2$. Set $k'=[k/2]$. For elements $b,b'$ in $B^{k,1}$ let 
$b\stackrel{\et{0}}{\longrightarrow}b'$ mean $\et{0}b=b'$. Then we have the following rules 
of $0$-actions. In (2)-(4) $\bullet$ stands for a crystal element whose explicit form is not used
later.

\noindent
(1)
\[
\begin{array}{c}
1\\ 2\\ \vdots\\ k
\end{array}
\stackrel{\et{0}}{\longrightarrow}
\begin{array}{c}
1\\ 3\\ \vdots\\ k\\ \ol{1}
\end{array}
\stackrel{\et{0}}{\longrightarrow}
\begin{array}{c}
3\\ 4\\ \vdots\\ k\\ \ol{2}\\ \ol{1}
\end{array}
\]
(2) Let $1\le m\le k'-1$.
\[
\begin{array}{c}
1\\ 2\\ \vdots\\ k-m\\ \ol{k-m}\\ \vdots\\ \ol{k-2m+1}
\end{array}
\stackrel{\et{0}}{\longrightarrow}
\begin{array}{c}
1\\ 2\\ \vdots\\ k-m-1\\ \ol{k-m-1}\\ \vdots\\ \ol{k-2m+1}\\ \ol{2}\\ \ol{1}
\end{array}
\stackrel{\et{0}}{\longrightarrow}
\hspace{5mm}\bullet
\]
(3) Let $0\le m_1\le m_2\le k'-1,m_1\le p\le \min(m_{21},m_s-1)$, where $m_{21}=m_2-m_1,m_s=m_1+m_2$.
\[
\begin{array}{c}
1\\ 2\\ \vdots\\ k-m_{21}-p\\ k-2p+1\\ \vdots\\ k-2m_1\\ \ol{k-m_{21}-p}\\ \vdots\\ \ol{k-2m_2+1}
\end{array}
\stackrel{\et{0}}{\longrightarrow}
\begin{array}{c}
1\\ 2\\ \vdots\\ k-m_{21}-p-1\\ k-2p+1\\ \vdots\\ k-2m_1\\ \ol{k-m_{21}-p-1}\\ \vdots\\ \ol{k-2m_2+1}\\
\ol{2}\\ \ol{1}
\end{array}
\stackrel{\et{0}}{\longrightarrow}
\hspace{5mm}\bullet
\]
(4) Let $1\le m_1\le m_2\le k'-1,m_{21}+1\le p\le m_2$, where $m_{21}=m_2-m_1$. Set $m_s=m_1+m_2$.
\[
\begin{array}{c}
1\\ 2\\ \vdots\\ k-m_s+p\\ \ol{k-m_s+p}\\ \vdots\\ \ol{k-2m_1+1}\\ \ol{k-2p}\\ \vdots\\ 
\ol{k-2m_2+1}
\end{array}
\stackrel{\et{0}}{\longrightarrow}
\begin{array}{c}
1\\ 2\\ \vdots\\ k-m_s+p-1\\ \ol{k-m_s+p-1}\\ \vdots\\ \ol{k-2m_1+1}\\ \ol{k-2p}\\ \vdots\\ 
\ol{k-2m_2+1}\\ \ol{2}\\ \ol{1}
\end{array}
\stackrel{\et{0}}{\longrightarrow}
\hspace{5mm}\bullet
\]
\end{lemma}

\begin{lemma} \label{lem:hwv}
Let $V$ be a $U_q(\geh_0)$-module with a crystal base $(L,B)$.
Let $W$ a submodule of $V$. Let $\{b_j\mid j\in I\}\subset B$. Suppose $v$ is a vector 
in $W$ such that $v\equiv\sum_j b_j\mbox{ mod }qL$. Decompose $I$ as $I=I_1\sqcup I_2$ by
\[
I_1=\{j\in I\mid \et{i}b_j=0\mbox{ for any }i\},\quad 
I_2=\{j\in I\mid \et{i}b_j\neq0\mbox{ for some }i\}.
\]
Then there exits a highest weight vector $w$ in $W$ such that $w\equiv\sum_{j\in I_1}b_j
\mbox{ mod }qL$.
\end{lemma}
\begin{proof}
By applying $\ft{i}\et{i}$ we know that there exists a vector $v'$ in $W$ such that 
$v'\equiv\sum_{j\in I'}b_j\mbox{ mod }qL$, where $I'=\{j\in I_2\mid\et{i}b_j\neq0\}$.
Hence there also exists a vector $v''$ in $W$ such that $v''\equiv\sum_{j\in (I')^c}b_j
\mbox{ mod }qL$. Continuing this with different $i$'s, we obtain a vector $v'''$ in $W$ 
such that $v'''\equiv\sum_{j\in I_1}b_j\mbox{ mod }qL$. Hence we can write $v'''$ as $v'''=w+w'$ 
in such a way that $w\equiv\sum_{j\in I_1}b_j\mbox{ mod }qL$ is a highest weight vector and
$w'\in qL$ is not, but we can remove $w'$ from $v'''$.
\end{proof}

In what follows in this section, by abuse of notation we represent a basis vector $v_I$ in 
$W_1^{(k)}$ also as $I$.

\begin{lemma}
Let $0\le m\le k'$. A highest weight vector of $V(2\olLa_{k-2m})$ in $W$ is given by
\[
\begin{array}{c}
1\\ 2\\ \vdots\\ k-m\\ \ol{k-m}\\ \vdots\\ \ol{k-2m+1}
\end{array}\ot
\begin{array}{c}
1\\ 2\\ \vdots\\ k-m\\ \ol{k-m}\\ \vdots\\ \ol{k-2m+1}
\end{array}
\quad\mbox{mod }qL.
\]
\end{lemma}

\begin{proof}
A highest weight vector of $V(2\olLa_k)$ is given by
\[
\begin{array}{c}
1\\ 2\\ \vdots\\ k
\end{array}
\ot
\begin{array}{c}
1\\ 2\\ \vdots\\ k
\end{array}.
\]
Noting that $W$ is a submodule of $W_{1,q^{-1}}^{(k)}\ot W_{1,q}^{(k)}$, apply $e_0^2$ and use 
Lemma \ref{lem:mod qL} and \ref{lem:e0} (1). We obtain 
\[
q
\begin{array}{c}
3\\ 4\\ \vdots\\ k\\ \ol{2}\\ \ol{1}
\end{array}
\ot
\begin{array}{c}
1\\ 2\\ \vdots\\ k
\end{array}
+(1+q^2)
\begin{array}{c}
1\\ 3\\ \vdots\\ k\\ \ol{1}
\end{array}
\ot
\begin{array}{c}
1\\ 3\\ \vdots\\ k\\ \ol{1}
\end{array}
+q
\begin{array}{c}
1\\ 2\\ \vdots\\ k
\end{array}
\ot
\begin{array}{c}
3\\ 4\\ \vdots\\ k\\ \ol{2}\\ \ol{1}
\end{array}
\equiv 
\begin{array}{c}
1\\ 3\\ \vdots\\ k\\ \ol{1}
\end{array}
\ot
\begin{array}{c}
1\\ 3\\ \vdots\\ k\\ \ol{1}
\end{array}
\quad\mbox{mod }qL
\]
as a vector in $W$. Apply further $e_{k-2}^{(2)}\cd e_2^{(2)}e_1^{(2)}
e_{k-1}^{(2)}\cd e_3^{(2)}e_2^{(2)}$, one obtains 
\[
\begin{array}{c}
1\\ 2\\ \vdots\\ k-1\\ \ol{k-1}
\end{array}
\ot
\begin{array}{c}
1\\ 2\\ \vdots\\ k-1\\ \ol{k-1}
\end{array}
\quad\mbox{mod }qL
\]
as a vector in $W$.

For $m>1$, we prove by induction on $m$. By Lemma \ref{lem:mod qL} and \ref{lem:e0} (2), one has
\begin{align*}
&e_{k-2m-2}^{(2)}\cd e_2^{(2)}e_1^{(2)}e_{k-2m-1}^{(2)}\cd e_3^{(2)}e_2^{(2)}e_0^2
\begin{array}{c}
1\\ 2\\ \vdots\\ k-m\\ \ol{k-m}\\ \vdots\\ \ol{k-2m+1}
\end{array}
\ot
\begin{array}{c}
1\\ 2\\ \vdots\\ k-m\\ \ol{k-m}\\ \vdots\\ \ol{k-2m+1}
\end{array}\\
&\equiv
\begin{array}{c}
1\\ 2\\ \vdots\\ k-m-1\\ \ol{k-m-1}\\ \vdots\\ \ol{k-2m-1}
\end{array}
\ot
\begin{array}{c}
1\\ 2\\ \vdots\\ k-m-1\\ \ol{k-m-1}\\ \vdots\\ \ol{k-2m-1}
\end{array}
\mbox{ mod }qL,
\end{align*}
as a vector in $W$. Lemma \ref{lem:hwv} completes the proof.
\end{proof}

\begin{lemma}
Let $0\le m_1\le m_2\le k'$. Set $m_{21}=m_2-m_1,m_s=m_1+m_2,M=\max(m_1,m_{21})$.
A highest weight vector of $V(\olLa_{k-2m_1}+\olLa_{k-2m_2})$ in $W$ is given by
\begin{equation} \label{hwv}
\sum_{p=m_1}^M
\begin{array}{c}
1\\ 2\\ \vdots\\ k-p\\ \ol{k-p}\\ \vdots\\ \ol{k-2p+1}
\end{array}
\ot
\begin{array}{c}
1\\ 2\\ \vdots\\ k-m_{21}-p\\ k-2p+1\\ \vdots\\ k-2m_1\\ \ol{k-m_{21}-p}\\ \vdots\\ \ol{k-2m_2+1}
\end{array}
+\sum_{p=M+1}^{m_2}
\begin{array}{c}
1\\ 2\\ \vdots\\ k-p\\ \ol{k-p}\\ \vdots\\ \ol{k-2p+1}
\end{array}
\ot
\begin{array}{c}
1\\ 2\\ \vdots\\ k-m_s+p\\ \ol{k-m_s+p}\\ \vdots\\ \ol{k-2m_1+1}\\ \ol{k-2p}\\ \vdots\\
\ol{k-2m_2+1}
\end{array}
\end{equation}
mod $qL$.
\end{lemma}

\begin{proof}
We prove by induction on $m_2$. The case of $m_2=m_1$ is proved in the previous lemma. Assume $m_1>0$.
Apply $e_{k-2m_2-2}\cd e_1e_{k-2m_2-1}\cd e_2e_0$ to \eqref{hwv} and use Lemma \ref{lem:mod qL} and 
\ref{lem:e0} (2),(3),(4). Since one can always neglect terms corresponding to the crystal elements
that are not killed by $\et{i}$ for some $i\neq0$ by Lemma \ref{lem:hwv}, we obtain
\begin{align*}
&\sum_{p=m_1}^M
\begin{array}{c}
1\\ 2\\ \vdots\\ k-p\\ \ol{k-p}\\ \vdots\\ \ol{k-2p+1}
\end{array}
\ot
\begin{array}{c}
1\\ 2\\ \vdots\\ k-m_{21}-p-1\\ k-2p+1\\ \vdots\\ k-2m_1\\ \ol{k-m_{21}-p-1}\\ \vdots\\ \ol{k-2m_2-1}
\end{array}
+\sum_{p=M+1}^{m_2}
\begin{array}{c}
1\\ 2\\ \vdots\\ k-p\\ \ol{k-p}\\ \vdots\\ \ol{k-2p+1}
\end{array}
\ot
\begin{array}{c}
1\\ 2\\ \vdots\\ k-m_s+p-1\\ \ol{k-m_s+p-1}\\ \vdots\\ \ol{k-2m_1+1}\\ \ol{k-2p}\\ \vdots\\
\ol{k-2m_2-1}
\end{array}\\
&+
\begin{array}{c}
1\\ 2\\ \vdots\\ k-m_2-1\\ \ol{k-m_2-1}\\ \vdots\\ \ol{k-2m_2-1}
\end{array}
\ot
\begin{array}{c}
1\\ 2\\ \vdots\\ k-m_1\\ \ol{k-m_1}\\ \vdots\\ \ol{k-2m_1+1}
\end{array}
\end{align*}
mod $qL$ as a vector in $W$. Note that the last term can be regarded as the term in the summand
of the middle term for $p=m_2+1$. If $m_1>m_{21}$, the induction proceeds. If $m_1\le m_{21}$, 
note also that the term in the summand of the middle term for $p=M+1$ can be regarded as the one 
of the first term for $p=M+1$.

The proof for $m_1=0$ is similar, but needs some attention. Note that only the first summation
survives in \eqref{hwv}. Divide the cases into three: $p=0,1\le p\le m_2-1,p=m_2$ for the 
calculation of the action of $e_0$.
\end{proof}

The following lemma is an easy consequence of \eqref{wRw'} with $l=2$ and the nondegeneracy of 
the admissible pairing.
\begin{lemma} \label{lem:characterizeN}
Let $(\;,\;)$ be the admissible pairing between $W_{1,q}^{(k)}\ot W_{1,q^{-1}}^{(k)}$ and
$W_{1,q^{-1}}^{(k)}\ot W_{1,q}^{(k)}$ induced from the poralization of $W_1^{(k)}$. Then
$u\in N$ if and only if $(u,v)=0$ for any $v\in W$.
\end{lemma}

Now we are to prove Proposition \ref{prop:N}. Set 
\[
b(\olLa_{k-2m})=12\cd(k-m)(\ol{k-m})\cd(\ol{k-2m+1})\in B.
\]
$b(\olLa_{k-2m})$ is a highest weight vector of $V(\olLa_{k-2m})$ in $W_1^{(k)}$.
We define the following subsets of $B^{\ot2}$.
\begin{align}
B_1^h&=\{b(\olLa_{k-2m_2})\ot b(\olLa_{k-2m_1})\mid 0\le m_1\le m_2\le k'\}, \label{B1h}\\
B_2^h&=\{\text{elements appearing in the summand of \eqref{hwv}}\}\setminus B_1^h, \label{B2h}\\
B_a&=\left(\bigcup_{i_1,\ldots,i_m\in I_0}\ft{i_m}\cd\ft{i_1}B_a^h\right)\setminus\{0\}
\quad\text{for }a=1,2, \label{Ba}\\
B_3&=B^{\ot2}\setminus(B_1\sqcup B_2),\quad 
B_3^h=\{b\in B_3\mid \et{i}b=0\text{ for any }i\neq0\}.
\end{align}
\noindent
Note that \eqref{hwv} always contains $b(\olLa_{k-2m_2})\ot b(\olLa_{k-2m_1})$.
For $J=I_1\ot I_2\in B^{\ot2}$ set $v_J=v_{I_1}\ot v_{I_2}$.

\begin{lemma}
For any $J\in B_2^h\sqcup B_3^h$ there exists a highest weight vector $w_J$ in $N$ such that 
\[
w_J\equiv\left\{
\begin{array}{ll}
v_J-v_{J_1}\quad&\text{if }J\in B_2^h\\
v_J&\text{if }J\in B_3^h
\end{array}\right.
\]
mod $qL^{\ot2}$. Here $J_1$ is the unique element of $B_1^h$ that has the same weight as $J$.
\end{lemma}
\begin{proof}
If the weight of $J$ is not of the form of $\olLa_{k-2m_1}+\olLa_{k-2m_2}$, the assertion is clear.
Suppose $\wt J=\olLa_{k-2m_1}+\olLa_{k-2m_2}$. Consider a highest weight vector of the form
\[
v=c_{J_1}v_{J_1}+\sum_{J_2}c_{J_2}v_{J_2}+\sum_{J_3}c_{J_3}v_{J_3}+v'.
\]
Here $J_1$ is the unique element in $B_1^h$ of the fixed weight, $c_{J_a}\in A$ for $a=1,2,3$,
the summation $\sum_{J_a}$ ranges over $J_a\in B_a^h$ that has the fixed weight for $a=2,3$, and
$v'$ is some vector in $qL^{\ot2}$. Let $w$ be the highest weight vector \eqref{hwv} in $W$.
{}From Lemma \ref{lem:characterizeN} $v\in N$ if and only if $(v,w)=0$. On the other hand, we have
$(v,w)\equiv c_{J_1}+\sum_{J_2}c_{J_2}$ mod $qA$ since $(L^{\ot2},L^{\ot2})\subset A$. Thus we have
\[
c_{J_1}+\sum_{J_2}c_{J_2}\equiv0\text{ mod }qA.
\]
Hence for any $c_{J_2},c_{J_3}\in A$, there exists a highest weight vector $v$ in $N$ such that
$v\equiv\sum_{J_2}c_{J_2}(v_{J_2}-v_{J_1})+\sum_{J_3}c_{J_3}v_{J_3}$ mod $qL^{\ot2}$. The assertion 
follows from this by setting $c_J=1,c_{J'}=0$ for $J'\neq J$.
\end{proof}

\noindent
Applying $\ft{i}$'s to $w_J$ for all $J$ in $B_2^h\sqcup B_3^h$ we obtain all weight vectors of $N$. 
Any such weight vector should have the following form
\[
av_J+\sum_{J'\ne J,J'\in B_2\sqcup B_3}c_{J'}v_{J'}+\sum_{J''\in B_1}d_{J''}v_{J''},
\]
where $J\in B_2\sqcup B_3,a\in1+qA,c_{J'}\in qA,d_{J''}\in A$. By Gaussian elimination, we obtain 
the desired result.

\section{Proof of proposition \ref{prop:norm}} \label{sec:pr2}

We begin this section with an easy lemma.

\begin{lemma} \label{lem:wt crit}
Let $v$ be a weight vector in $W_l^{(k)}$. If $(\wt v,\epsilon_i)>l$ for some $i\in\{1,\ldots,n\}$, 
then $v=0$.
\end{lemma}
\begin{proof}
The claim follows from the fact that $(\wt u,\eps_i)\le1$ for a nonzero weight vector $u$ in 
$W_1^{(k)}$ and $W_l^{(k)}$ is a subspace of $W_1^{(k)}$.
\end{proof}

\noindent
The following formula will be used frequently.
\begin{equation} \label{PBW}
f_i^{(a)}e_i^{(b)}=\sum_{j=0}^{\min(a,b)}e_i^{(b-j)}f_i^{(a-j)}{q^{a-b}t_i^{-1}\brace j},
\end{equation}
where ${t\brace j}=\prod_{k=1}^j(q^{1-k}t-q^{k-1}t^{-1})/(q^k-q^{-k})$.

In this section we abbreviate $l$ of the preporlarization $(\;,\;)_l$ on $W_l^{(k)}$.
We also write $|u|^2$ for $(u,u)$. Recall the definition of $u_m$ \eqref{u_m}. 
$\wt u_m$ is given by
\[
\wt u_m=(l-c_1)\olLa_k+(c_1-c_2)\olLa_{k-2}+\cd+(c_{m-1}-c_m)\olLa_{k-2m+2}+c_m\olLa_{k-2m}.
\]

\begin{lemma} \label{lem:rec}
\[
|u_m|^2=q^{c_m(2l-c_m)}{2l\brack c_m}|u_{m-1}|^2.
\]
\end{lemma}

\begin{proof}
Since the other case is similar, we prove when $k$ is even. Using \eqref{admissible2}, we have
\[
|u_m|^2=((e_{k-2m-1}^{(c_m)}\cd e_1^{(c_m)})(e_{k-2m+1}^{(c_m)}\cd e_2^{(c_m)})e_0^{(c_m)}u_{m-1},
f_{k-2m}^{(c_m)}u_m).
\]
By \eqref{PBW} we obtain
\begin{equation} \label{pr1}
f_{k-2m}^{(c_m)}u_m=\sum_je_{k-2m}^{(c_m-j)}f_{k-2m}^{(c_m-j)}{c_m\brack j}
(e_{k-2m-1}^{(c_m)}\cd e_1^{(c_m)})(e_{k-2m+1}^{(c_m)}\cd e_2^{(c_m)})e_0^{(c_m)}u_{m-1}.
\end{equation}
Note that $\wt f_{k-2m}^{(c_m-j)}(e_{k-2m-1}^{(c_m)}\cd e_1^{(c_m)})
(e_{k-2m+1}^{(c_m)}\cd e_2^{(c_m)})e_0^{(c_m)}u_{m-1}=\wt u_m-(2c_m-j)\alpha_{k-2m}$. From Lemma
\ref{lem:wt crit} the summand in the r.h.s. of \eqref{pr1} becomes $0$ unless $j=c_m$. Hence we have
\[
|u_m|^2=|(e_{k-2m-1}^{(c_m)}\cd e_1^{(c_m)})(e_{k-2m+1}^{(c_m)}\cd e_2^{(c_m)})e_0^{(c_m)}u_{m-1}|^2.
\]
Similar calculations continue until we arrive at $|u_m|^2=|e_0^{(c_m)}u_{m-1}|^2$. Using 
\eqref{admissible2},\eqref{PBW} and Lemma \ref{lem:wt crit} again, we this time have
$|u_m|^2=q^{c_m(2l-c_m)}{2l\brack c_m}|u_{m-1}|^2$.
\end{proof}

\begin{lemma} \label{lem:e_j}
\begin{itemize}
\item[(1)] $e_ju_{k'}=0$ when $k$ is even, if $j\ge k+1$ or $j=1$ when $k$ is odd.
\item[(2)] $|f_ju_p|^2=q^{c_p(2l-1-c_p)}{2l-1\brack c_p}q^{c_{p-1}-1}[c_{p-1}]|u_p|^2$
	if $j=k-2p+2,p=1,\ldots,k'$.
\item[(3)] $|f_ju_p|^2=0$ if $j=k-2p+1,p=1,\ldots,k'$.
\end{itemize}
\end{lemma}

\begin{proof}
(1) Write $u_{k'}=Eu_0$. If $j\ge k+1$, $e_j$ commutes with $E$. The claim follows from $e_ju_0=0$.
When $k$ is odd, $e_1u_{k'}=0$ follows from Lemma \ref{lem:wt crit}.

(2) When $c_p=0$, the equality is shown as follows.
\begin{align*}
|f_ju_p|^2&=|f_{k-2p+2}u_{p-1}|^2\\
&=(u_{p-1},q^{-1}t_{k-2p+2}e_{k-2p+2}f_{k-2p+2}u_{p-1})\\
&=q^{c_{p-1}-1}[c_{p-1}]|u_{p-1}|^2.
\end{align*}
Here we have used the relation $e_{k-2p+2}u_{p-1}=0$, that can be confirmed by Lemma 
\ref{lem:wt crit}.

Now assume $c_p>0$. Imitating the proof of Lemma \ref{lem:rec}, one obtains 
\[
|f_ju_p|^2=|(e_{k-2p+1}^{(c_p)}\cd e_2^{(c_p)})e_0^{(c_p)}\cdot f_{f-2p+2}u_{p-1}|^2.
\]
Next we calculate 
\begin{align*}
&q^{-c_p^2}t_{k-2p+1}^{-c_p}f_{k-2p+1}^{(c_p)}(e_{k-2p+1}^{(c_p)}\cd e_2^{(c_p)})e_0^{(c_p)}\cdot
f_{k-2p+2}u_{p-1}\\
&=q^{-c_p}\sum_je_{k-2p+1}^{(c_p-j)}f_{k-2p+1}^{(c_p-j)}{c_p-1\brack j}
e_{k-2p}^{(c_p)}\cd e_2^{(c_p)}e_0^{(c_p)}\cdot f_{k-2p+2}u_{p-1}.
\end{align*}
Since $(\wt f_{k-2p+1}^{(c_p-j)}e_{k-2p}^{(c_p)}\cd e_2^{(c_p)}e_0^{(c_p)}f_{k-2p+2}u_{p-1},
\eps_{k-2p+2})=l-1+c_p-j$, the summand of the above expression survives only when $j=c_p,c_p-1$.
Noting ${c_p-1\brack c_p}=0$, we obtain
\[
|f_ju_p|^2=|(e_{k-2p}^{(c_p)}\cd e_2^{(c_p)})e_0^{(c_p)}\cdot f_{k-2p+1}f_{k-2p+2}u_{p-1}|^2.
\]
Calculating similarly, one gets $|f_ju_{p-1}|^2=|e_0^{(c_p)}\cdot f_2\cd f_{k-2p+2}u_{p-1}|^2$.
After removing $e_0^{(c_p)}$, we arrive at 
\[
|f_ju_p|^2=q^{c_p(2l-1-c_p)}{2l-1\brack c_p}|f_2\cd f_{k-2p+2}u_{p-1}|^2.
\]
Calculating similarly, we obtain
\[
|f_2\cd f_{k-2p+2}u_{p-1}|^2=|f_{k-2p+2}u_{p-1}|^2=q^{c_{p-1}-1}[c_{p-1}]|u_{p-1}|^2.
\]

(3) The proof goes parallel to that of (2). When $c_p=0$, 
\[
|f_ju_p|^2=|f_{k-2p+1}u_{p-1}|^2=0
\]
from Lemma \ref{lem:wt crit}. Assume $c_p>0$. One obtains 
\[
|f_ju_p|^2=|f_1f_2\cd f_{k-2p+1}u_{p-1/2}|^2.
\]
Noting that $f_iu_{p-1}=0$ for $i=1,\ldots,k-2p+1$, we have
\begin{align*}
f_1f_2\cd f_{k-2p+1}u_{p-1/2}&=f_1f_2\cd f_{k-2p+1}(e_{k-2p+1}^{(c_p)}\cd e_2^{(c_p)})
e_0^{(c_p)}u_{p-1}\\
&=\alpha\cdot f_1e_{k-2p+1}^{(c_p-1)}\cd e_2^{(c_p-1)})\cdot e_0^{(c_p)}u_{p-1}=0.
\end{align*}
Here $\alpha$ is a product of $q$-integers.

The proof is complete.
\end{proof}

Now we are in a position to prove Proposition \ref{prop:norm}.
(1) is a simple consequence of Lemma \ref{lem:rec}. (2) when $j\ge k+1$ is settled by Lemma
\ref{lem:e_j} (1). To show when $j\le k$ note that 
\begin{align*}
|e_ju_{k'}|^2&=(u_{k'},q^{-1}t_j^{-1}f_je_ju_{k'})\\
&=q^{2\beta_j}|f_ju_{k'}|^2+q^{\beta_j-1}[\beta_j]|u_{k'}|^2,
\end{align*}
where
\[
\beta_j=-\langle h_j,\wt u_{k'}\rangle=\left\{
\begin{array}{ll}
c_{\frac{k-j}2+1}-c_{\frac{k-j}2}\quad&\mbox{if }j\equiv k\:(2),\\
0&\mbox{if }j\not\equiv k\:(2).
\end{array}\right.
\]
Thus we are left to evaluate $|f_ju_{k'}|^2$. Examining the proof of Lemma \ref{lem:rec} carefully,
one notices that the same recursion formula is valid when $m>p$, namely, one has
\[
|f_ju_m|^2=q^{c_m(2l-c_m)}{2l\brack c_m}|f_ju_{m-1}|^2\quad\mbox{for }m>p.
\]
The formula for $|f_ju_{k'}|^2$ is obtained from this, Lemma \ref{lem:e_j} (2) or (3) and Lemma
\ref{lem:rec}. Calculating explicitly we obtain (2).

\section*{Acknowledgments}
\smallskip\par\noindent
The author thanks Masaki Kashiwara, Hiraku Nakajima, Satoshi Naito and Daisuke Sagaki for 
stimulating discussions.
He is partially supported by
Grant-in-Aid for Scientific Research (C) 18540030,
Japan Society for the Promotion of Science.

\end{document}